\documentclass[11pt]{amsart}
\usepackage{amssymb,amsmath,amsthm}
\usepackage[margin=1.15in]{geometry}
\usepackage{enumitem}
\usepackage{xcolor}
\usepackage[colorlinks=true,linkcolor=blue!50!black,citecolor=blue!50!black,urlcolor=blue!50!black]{hyperref}

\newtheorem{theorem}{Theorem}[section]
\newtheorem{lemma}[theorem]{Lemma}
\newtheorem{proposition}[theorem]{Proposition}
\newtheorem{corollary}[theorem]{Corollary}
\newtheorem{conjecture}[theorem]{Conjecture}
\theoremstyle{definition}
\newtheorem{definition}[theorem]{Definition}
\theoremstyle{remark}
\newtheorem{remark}[theorem]{Remark}

\numberwithin{equation}{section}

\newcommand{\N}{\mathbb{N}}
\newcommand{\Z}{\mathbb{Z}}
\newcommand{\Q}{\mathbb{Q}}
\newcommand{\R}{\mathbb{R}}
\newcommand{\eps}{\varepsilon}
\newcommand{\abs}[1]{\left|#1\right|}

\makeatletter
\@namedef{subjclassname@2020}{\textup{2020} Mathematics Subject Classification}
\makeatother
\date{Revised version, November, 2026}
\begin{document}

\title[Dirichlet approximation for affine actions]{A generalization of the Dirichlet approximation theorem\\ for affine actions on the real line}

\author{Mohammad Javaheri}
\address{Department of Mathematics and Statistics, Siena University, Loudonville, NY 12211}
\email{mjavaheri@siena.edu}


\begin{abstract}
We study quantitative density of orbits of the free semigroup $\Gamma$ generated by two maps of $(0,\infty)$. For the M\"obius pair $x\mapsto 1/x$, $x\mapsto x+1$ we show that every orbit is $O(n^{-1})$-dense in $[0,\infty)$, that every irrational number $u$ is approached at a rate governed by the continued fraction of $u$ (at least $n^{-2}$, exponential when $u$ has bounded partial quotients, and $\exp(-cn/\log n)$ for almost every $u$), and that for rational $u$ and rational starting points the rate $n^{-1}$ is optimal. For the affine pair $T_0(x)=x/a$, $T_1(x)=bx+1$ with $a,b>1$, which was studied by Bergelson, Misiurewicz and Senti, we prove that for every $x>0$ and $u\ge 0$ there are infinitely many $\gamma\in\Gamma$ with
\[
|\gamma(x)-u|<C\,\lambda^{-|\gamma|},\qquad \lambda=\exp\!\Big(\frac{\log a\,\log\frac{a}{a-1}}{\log\frac{a^2b}{a-1}}\Big)>1,
\]
where $|\gamma|$ is the word length of $\gamma$ and $C$ is independent of $\gamma$. This exponential rate holds whether or not $\log_a b$ is rational, and the exponential order cannot be improved for almost every $u$. We also give a corrected version of the Kronecker-type estimate $|\gamma(x)-u|<(u+1)(t^{1/|\gamma|}-1)$. Finally, we discuss the approximation of real numbers by ratios of primes and by ratios of logarithms of primes, updated in light of recent results on gaps between primes.
\end{abstract}

\maketitle

\section{Introduction}

The Dirichlet approximation theorem states that, given a real number $u$ and a positive integer $N$, there exist integers $p,q$ with $0<q\le N$ such that
\begin{equation}\label{eq:dirichlet}
\abs{\frac pq-u}<\frac{1}{qN}.
\end{equation}
In particular, if $u$ is irrational then $|p/q-u|<q^{-2}$ for infinitely many reduced fractions $p/q$. We give a formulation of this theorem in terms of semigroup actions which leads to natural generalizations.

Let $\Gamma$ be the free semigroup generated by two letters $\sigma$ and $\tau$, with identity $\lambda_0$ (the empty word). Given two maps $T_0,T_1:(0,\infty)\to(0,\infty)$, we obtain an action of $\Gamma$ on $(0,\infty)$ by letting $\sigma$ act as $T_0$ and $\tau$ act as $T_1$; words act as compositions, so that $(w_1w_2)(x)=w_1(w_2(x))$ (the rightmost letter acts first). We write $\gamma\in\Gamma(T_0,T_1)$ when we want to emphasize the choice of $T_0$ and $T_1$. The \emph{length} $|\cdot|:\Gamma\to\Z_{\ge0}$ is defined by $|\lambda_0|=0$, $|\sigma|=|\tau|=1$ and $|w_1w_2|=|w_1|+|w_2|$. For $x,y>0$ we say that $y$ is \emph{reachable} from $x$ if $y=\gamma(x)$ for some $\gamma\in\Gamma$, and we define the reachability function
\begin{equation}\label{eq:R}
R(x,y)=\min\{|\gamma| : \gamma\in\Gamma,\ \gamma(x)=y\},
\end{equation}
with $R(x,y)=\infty$ if $y$ is not reachable from $x$.

\begin{definition}\label{def:dense}
For $x>0$ let $\Gamma x$ denote the orbit of $x$ under $\Gamma(T_0,T_1)$. Let $f:\N\to(0,\infty)$ be a non-increasing function with $f(n)\to0$. We say that $\Gamma x$ is \emph{$O(f)$-dense} in $[0,\infty)$ if for every $u\ge0$ there exist a constant $C=C(T_0,T_1,x,u)>0$ and infinitely many $\gamma\in\Gamma$ such that
\[
|\gamma(x)-u|<C f(|\gamma|).
\]
\end{definition}

Since only finitely many words have a given length, ``infinitely many $\gamma$'' forces $|\gamma|\to\infty$ along the approximating words.

\subsection*{The M\"obius pair}
In Section~\ref{sec:mobius} we take
\begin{equation}\label{eq:mobius}
T_0(x)=\frac1x,\qquad T_1(x)=x+1.
\end{equation}
Here $\gamma(1)$ is always a positive rational number, and Lemma~\ref{lem:R} shows $R(1,p/q)\le p+q$; combined with \eqref{eq:dirichlet} this shows that every irrational $u$ is approached at rate $O(|\gamma|^{-2})$. We prove (Theorem~\ref{thm:mobius}) that every orbit is $O(n^{-1})$-dense in $[0,\infty)$, with the rate $n^{-2}$ at irrational targets. We then refine this in two directions. First (Theorem~\ref{thm:cf} and Corollary~\ref{cor:cf}), the true rate at an irrational $u$ is governed by the continued fraction expansion of $u$: it is exponential in $|\gamma|$ when $u$ has bounded partial quotients, and it is at least $\exp(-c|\gamma|/\log|\gamma|)$ for Lebesgue almost every $u$, for every $c<\pi^2/6$. Second (Proposition~\ref{prop:sharp}), for rational targets and rational starting points the exponent $1$ is optimal: by an argument with the Stern--Brocot tree, $|\gamma(x)-k/l|\ge 1/(l^2(s(x)+|\gamma|))$ whenever $\gamma(x)\ne k/l$, where $s(x)$ is the sum of the partial quotients of $x$.

\subsection*{The affine pair}
In Section~\ref{sec:affine} we consider
\begin{equation}\label{eq:affine}
T_0(x)=\frac xa,\qquad T_1(x)=bx+1,\qquad a,b>1.
\end{equation}
Bergelson, Misiurewicz and Senti \cite{BMS} showed that every orbit of this action is dense in $[0,\infty)$. Our main result (Theorem~\ref{thm:exp}) is that every orbit is in fact \emph{exponentially} dense: for every $x>0$ and $u\ge0$ there are infinitely many $\gamma$ with $|\gamma(x)-u|<C\lambda^{-|\gamma|}$ where $\lambda=\lambda(a,b)>1$ is explicit. The proof rests on an explicit normal form for $\gamma(x)$ (Lemma~\ref{lem:normal}) and a greedy algorithm, and it does not distinguish between rational and irrational values of $\log_ab$. A counting argument (Proposition~\ref{prop:count}) shows that for almost every $u$ no rate faster than $\mu^{-|\gamma|}$ with $\mu>2$ is possible, so the exponential order is optimal.

We also give corrected versions of two earlier estimates: a Kronecker-type bound $|\gamma(x)-u|<(u+1)(t^{1/|\gamma|}-1)$ with $t=a^{12}b^4$ when $\log_ab$ is irrational (Theorem~\ref{thm:kronecker}), and an exponential bound with exponent independent of $u$ when $\log_ab$ is rational (Theorem~\ref{thm:rational}). Since $t^{1/n}-1\sim(\log t)/n$, the Kronecker-type bound is an $O(n^{-1})$ statement; it is superseded by Theorem~\ref{thm:exp}, but its proof is short and the words it produces are very simple.

The contrast between the two sections is worth noting: for the M\"obius pair rational targets are approached only at rate $n^{-1}$ (at least from rational starting points), whereas for the affine pair every target is approached exponentially fast.

\subsection*{Relation to the $3x+1$ problem}
The $3x+1$ (Collatz) problem asks whether iterating
\[
f(m)=\begin{cases} m/2, & m \text{ even},\\ 3m+1, & m\text{ odd},\end{cases}
\]
eventually reaches $1$ from every positive integer; see \cite{Lag} for a survey. In our notation, with $a=2$ and $b=3$, a Collatz trajectory is a particular orbit segment of $\Gamma(T_0,T_1)$, namely one in which $T_0$ is applied exactly to even integers and $T_1$ exactly to odd integers. Thus the Collatz conjecture implies that $1$ is reachable from every positive integer via $\Gamma(T_0,T_1)$; the converse is not clear, since the semigroup allows $T_1$ to be applied to even integers as well.

Moreover, any sequence $x_1,\dots,x_m$ with $x_1\in\N$, $x_m=1$ and $x_{i+1}\in\{T_0(x_i),T_1(x_i)\}$ consists of integers. Indeed, write $v_2$ for the $2$-adic valuation on $\Q^\times$. If $v_2(x_i)<0$ then $v_2(T_0x_i)=v_2(x_i)-1<0$ and $v_2(T_1x_i)=v_2(3x_i+1)=v_2(x_i)<0$, so once a non-integer appears, all later terms are non-integers; since $x_m=1$ is an integer, no $x_i$ can be a non-integer.

Theorem~\ref{thm:exp} shows that, starting from \emph{any} $x>0$, there are words $\gamma$ of arbitrarily large length $n$ with $|\gamma(x)-1|<C\lambda^{-n}$, where $\lambda=\lambda(2,3)\approx 1.2133$.

\subsection*{Primes}
In Section~\ref{sec:primes} we discuss approximating real numbers by ratios of primes and by ratios of logarithms of primes, updating the discussion in light of the work of Goldston, Pintz and Y\i ld\i r\i m \cite{GPY,GMPY}, Zhang \cite{Zhang}, Maynard \cite{Maynard} and Baker, Harman and Pintz \cite{BHP}.

\section{The M\"obius pair}\label{sec:mobius}

Throughout this section $T_0(x)=1/x$ and $T_1(x)=x+1$.

\subsection{Matrices and continued fractions}
Associate to $T_0$ and $T_1$ the matrices
\[
M_\sigma=\begin{pmatrix}0&1\\1&0\end{pmatrix},\qquad M_\tau=\begin{pmatrix}1&1\\0&1\end{pmatrix},
\]
and to a word $\gamma=w_1\cdots w_n$ the product $M_\gamma=M_{w_1}\cdots M_{w_n}$. If $M_\gamma=\begin{pmatrix}A&B\\C&D\end{pmatrix}$ then $A,B,C,D$ are non-negative integers, $AD-BC=\pm1$, $C+D\ge1$, and
\[
\gamma(z)=\frac{Az+B}{Cz+D}\qquad(z>0).
\]

Every positive rational number $r$ has a unique \emph{canonical} continued fraction expansion $r=[a_0;a_1,\dots,a_k]$ with $a_0\ge0$, $a_1,\dots,a_k\ge1$, and either $k=0$ (and then $a_0\ge1$) or $k\ge1$ and $a_k\ge2$. We put
\[
s(r)=a_0+a_1+\cdots+a_k,\qquad \ell(r)=k.
\]

\begin{lemma}\label{lem:s}
For every positive rational $r$ we have $s(T_0r)=s(r)$ and $s(T_1r)=s(r)+1$. Consequently, for every positive rational $x$ and every $\gamma\in\Gamma$,
\[
s(\gamma(x))=s(x)+(\text{number of occurrences of }\tau\text{ in }\gamma)\le s(x)+|\gamma|.
\]
\end{lemma}

\begin{proof}
Let $r=[a_0;a_1,\dots,a_k]$ be canonical. Then $T_1r=[a_0+1;a_1,\dots,a_k]$ is canonical. If $a_0=0$ then $k\ge1$ and $T_0r=[a_1;a_2,\dots,a_k]$ is canonical. If $a_0\ge1$ and $k\ge1$, then $T_0r=[0;a_0,a_1,\dots,a_k]$ is canonical; if $a_0\ge2$ and $k=0$ then $T_0r=[0;a_0]$ is canonical; and if $r=[1]$ then $T_0r=r$. In all cases the sum of the partial quotients is preserved by $T_0$ and increased by one by $T_1$.
\end{proof}

\begin{lemma}\label{lem:R}
Let $r=[a_0;a_1,\dots,a_k]$ be a positive rational number in canonical form. Then
\begin{equation}\label{eq:Rcf}
s(r)-1\le R(1,r)\le s(r)+\ell(r)-1 .
\end{equation}
In particular, if $p,q$ are positive integers then $R(1,p/q)\le p+q$.
\end{lemma}

\begin{proof}
The lower bound follows from Lemma~\ref{lem:s}, since $s(1)=1$ and each word reaching $r$ from $1$ contains at least $s(r)-1$ letters $\tau$. For the upper bound, put $z_k=a_k$ and $z_{i}=a_i+1/z_{i+1}$ for $i<k$, so that $z_0=r$. Then $z_k=T_1^{a_k-1}(1)$ and $z_i=T_1^{a_i}T_0(z_{i+1})$, so $r$ is reached from $1$ by a word of length $(a_k-1)+\sum_{i<k}(a_i+1)=s(r)+k-1$.

For the last statement we argue by induction on $p+q$ (one may also deduce it from \eqref{eq:Rcf}). If $p=q$ then $R(1,p/q)=0$. If $p>q$ then $p/q=T_1((p-q)/q)$, so $R(1,p/q)\le 1+R(1,(p-q)/q)\le 1+p$. If $1<p<q$ then $p/q=T_0T_1((q-p)/p)$, so $R(1,p/q)\le 2+R(1,(q-p)/p)\le q+2\le p+q$. Finally, if $p=1$ then $1/q=T_0T_1^{q-1}(1)$ and $R(1,1/q)\le q$.
\end{proof}

\begin{lemma}\label{lem:contract}
Let $\gamma\in\Gamma$ and write $\gamma(1)=p/q$ in lowest terms. Then for every $\eps\ge0$,
\[
|\gamma(1+\eps)-\gamma(1)|=\frac{\eps}{q(q+C\eps)}\le\frac{\eps}{q^2},
\]
where $C$ is the lower left entry of $M_\gamma$.
\end{lemma}

\begin{proof}
With $M_\gamma=\begin{pmatrix}A&B\\C&D\end{pmatrix}$ we have $\gamma(1)=(A+B)/(C+D)$, and this fraction is in lowest terms because $AD-BC=\pm1$; thus $p=A+B$ and $q=C+D$. Then
\[
\gamma(1+\eps)-\gamma(1)=\frac{(A+B)+A\eps}{(C+D)+C\eps}-\frac{A+B}{C+D}=\frac{\eps(AD-BC)}{(C+D)(C+D+C\eps)},
\]
which gives the claim.
\end{proof}

\subsection{Density}
Define $Q:[0,\infty)\to\{1,2\}$ by $Q(u)=1$ if $u$ is rational and $Q(u)=2$ if $u$ is irrational.

\begin{theorem}\label{thm:mobius}
For every $x>0$ the orbit $\Gamma x$ is $O(n^{-1})$-dense in $[0,\infty)$. More precisely, for every $x>0$ and $u\ge0$ there exist $C=C(x,u)>0$ and infinitely many $\gamma\in\Gamma$ such that
\begin{equation}\label{eq:mobius-rate}
|\gamma(x)-u|<C|\gamma|^{-Q(u)}.
\end{equation}
\end{theorem}

\begin{proof}
\emph{Reduction.} If $x<1$, replace $x$ by $T_1(x)=x+1$; a word $\gamma$ for $x+1$ gives the word $\gamma\tau$ for $x$, of length $|\gamma|+1\le2|\gamma|$. So we may assume $x=1+\eps$ with $\eps\ge0$.

\emph{The case $u=0$.} For $m\ge1$ let $\gamma=\sigma\tau^m$. Then $\gamma(x)=1/(x+m)<1/m<2/|\gamma|$.

\emph{Rational $u=k/l>0$.} For $m\ge1$ put $r_m=(km+1)/(lm)$, so $r_m-u=1/(lm)$. By Lemma~\ref{lem:R} there is $\gamma_m$ with $\gamma_m(1)=r_m$ and $|\gamma_m|\le(k+l)m+1\le(k+l+1)m$. Since $\gcd(km+1,m)=1$, the reduced denominator $q_m$ of $r_m$ is at least $lm/\gcd(km+1,l)\ge m$. By Lemma~\ref{lem:contract},
\[
|\gamma_m(x)-u|\le\frac{\eps}{q_m^2}+\frac1{lm}\le\frac{\eps+1}{m}\le\frac{x(k+l+1)}{|\gamma_m|}.
\]
The $r_m$ are distinct, so the $\gamma_m$ are distinct. (In the original version of this argument one took $r_m=(kl^m+1)/l^{m+1}$; this fails when $u$ is an integer, since then $l=1$ and $r_m$ does not depend on $m$.)

\emph{Irrational $u$.} This follows from Corollary~\ref{cor:cf}(a) below.
\end{proof}

\subsection{Irrational targets and continued fractions}
Let $u=[a_0;a_1,a_2,\dots]$ be irrational, with convergents $p_k/q_k$, and put $S_k=a_0+a_1+\cdots+a_k$.

\begin{theorem}\label{thm:cf}
Let $u>0$ be irrational and $x>0$. For every $k\ge1$ there is $\gamma_k\in\Gamma$ with
\[
|\gamma_k|\le S_k+k\qquad\text{and}\qquad |\gamma_k(x)-u|<\frac{x+1}{q_k^2}.
\]
\end{theorem}

\begin{proof}
The canonical expansion of $p_k/q_k$ is either $[a_0;\dots,a_k]$ (if $a_k\ge2$) or $[a_0;\dots,a_{k-1}+1]$ (if $a_k=1$); in both cases its partial quotients sum to $S_k$ and it has at most $k$ terms after $a_0$. By Lemma~\ref{lem:R} there is $\delta_k$ with $\delta_k(1)=p_k/q_k$ and $|\delta_k|\le S_k+k-1$. Recall $|p_k/q_k-u|<1/(q_kq_{k+1})\le 1/q_k^2$.

If $x=1+\eps\ge1$, let $\gamma_k=\delta_k$; by Lemma~\ref{lem:contract}, $|\gamma_k(x)-u|<(\eps+1)/q_k^2=x/q_k^2$. If $x<1$, let $\gamma_k=\delta_k\tau$, so that $\gamma_k(x)=\delta_k(1+x)$ and $|\gamma_k(x)-u|<(1+x)/q_k^2$.
\end{proof}

\begin{corollary}\label{cor:cf}
Let $x>0$ and let $u=[a_0;a_1,\dots]$ be irrational.
\begin{enumerate}[label=\textup{(\alph*)}]
\item There is $C=C(x,u)$ such that $|\gamma(x)-u|<C|\gamma|^{-2}$ for infinitely many $\gamma$.
\item If $a_i\le A$ for all $i\ge1$, then $|\gamma(x)-u|<C\varphi^{-2|\gamma|/(A+1)}$ for infinitely many $\gamma$, where $\varphi=(1+\sqrt5)/2$ and $C=C(x,u)$.
\item For Lebesgue almost every $u>0$ and every $c<\pi^2/6$, there are infinitely many $\gamma$ with $|\gamma(x)-u|<\exp(-c|\gamma|/\log|\gamma|)$.
\end{enumerate}
\end{corollary}

\begin{proof}
(a) We have $q_1=a_1$ and, for $k\ge2$, $q_k=a_kq_{k-1}+q_{k-2}\ge (a_k+q_{k-1}-1)+q_{k-2}\ge a_k+q_{k-1}$, using $(a_k-1)(q_{k-1}-1)\ge0$ and $q_{k-2}\ge1$. Hence $q_k\ge a_1+\cdots+a_k=S_k-a_0\ge k$, and
\[
|\gamma_k|\le S_k+k\le a_0+2q_k.
\]
By Theorem~\ref{thm:cf}, $|\gamma_k(x)-u|<(x+1)q_k^{-2}\le 4(x+1)(|\gamma_k|-a_0)^{-2}$, which is at most $C|\gamma_k|^{-2}$ for $k$ large. The words $\gamma_k$ are pairwise distinct, because the values $\delta_k(1)=p_k/q_k$ are distinct.

(b) Now $q_k\ge\varphi^{k-1}$ and $|\gamma_k|\le a_0+(A+1)k$, so $q_k^{-2}\le \varphi^{2}\varphi^{-2k}\le \varphi^{2+2a_0/(A+1)}\varphi^{-2|\gamma_k|/(A+1)}$.

(c) By L\'evy's theorem, $\frac1k\log q_k\to\beta:=\pi^2/(12\log2)$ for almost every $u$. By Khinchin's weak law for the sums of partial quotients, $S_k/(k\log k)\to 1/\log2$ in measure; hence there is a fixed sequence $k_j\to\infty$ along which $S_{k_j}/(k_j\log k_j)\to1/\log2$ almost everywhere (see \cite{IK,Khin}). Fix $\eta>0$ and such a $u$. For $j$ large, $n_j:=|\gamma_{k_j}|\le(1+\eta)k_j\log k_j/\log2$, hence $k_j\ge \frac{\log 2}{1+\eta}\,\frac{n_j}{\log n_j}$ (as $k_j\le n_j$), and
\[
|\gamma_{k_j}(x)-u|<(x+1)e^{-2(1-\eta)\beta k_j}\le (x+1)\exp\Big(-\frac{\pi^2}{6}\,\frac{1-\eta}{1+\eta}\,\frac{n_j}{\log n_j}\Big).
\]
Choosing $\eta$ small in terms of $c<\pi^2/6$ gives the claim.
\end{proof}

\begin{remark}
Part (a) is the estimate $|\gamma(x)-u|\le C|\gamma|^{-2}$ asserted for irrational $u$ in Theorem~\ref{thm:mobius}. Parts (b) and (c) show that it is far from optimal for most $u$; for instance, for the golden ratio the orbit of $1$ approaches $\varphi$ exponentially fast. On the other hand, by Proposition~\ref{prop:count} below, for almost every $u$ no rate $\mu^{-|\gamma|}$ with $\mu>2$ is possible.
\end{remark}

\subsection{Optimality for rational targets}
We use the Stern--Brocot tree (see \cite[\S4.5]{GKP}). Its nodes are the positive rationals, each occurring exactly once; the root is $1/1$, and the node $r=[a_0;a_1,\dots,a_k]$ (canonical) has depth $s(r)-1$. Every node $u=k/l$ has a nearest left ancestor $k_L/l_L<u$ and a nearest right ancestor $k_R/l_R>u$ (with the conventions $0/1$ and $1/0$ at the boundary) such that $k=k_L+k_R$, $l=l_L+l_R$ and $kl_L-lk_L=k_Rl-l_Rk=1$; moreover, the rationals in the open interval $(k_L/l_L,k_R/l_R)$ are exactly the nodes of the subtree rooted at $u$.

\begin{lemma}\label{lem:SB}
Let $u=k/l\ge0$ be rational in lowest terms and let $r>0$ be rational with $r\ne u$. Then
\[
|r-u|\ge\frac{1}{l^2\,s(r)}.
\]
\end{lemma}

\begin{proof}
If $u=0$ (so $l=1$): if $r\ge1$ then $r\ge1\ge 1/s(r)$; otherwise $r=[0;a_1,\dots,a_m]$ with $s(r)=a_1+\dots+a_m$, and $r=1/a_1\ge 1/s(r)$ if $m=1$, while $r>1/(a_1+1)\ge 1/s(r)$ if $m\ge2$ (as then $a_m\ge 2$).

Let $u>0$, let $d_0$ be the depth of $u$, and let $D=s(r)-1$ be the depth of $r$. If $r\notin(k_L/l_L,k_R/l_R)$ then $|r-u|\ge\min\{1/(ll_L),1/(ll_R)\}\ge1/l^2$ since $l_L,l_R\le l$. Otherwise $r$ lies in the subtree of $u$; suppose $r<u$ (the case $r>u$ is symmetric). The nodes $c_j=(k_L+jk)/(l_L+jl)$, $j\ge1$, form the right spine of the left subtree of $u$; $c_j$ has depth $d_0+j$, and the rationals in $(c_j,u)$ are exactly the nodes of the subtree rooted at $c_{j+1}$, all of depth at least $d_0+j+1$. Here $D\ge d_0+1$; taking $j=D-d_0\ge1$, the node $r$ is not in $(c_j,u)$, so $r\le c_j$ and therefore
\[
u-r\ge u-c_j=\frac{kl_L-lk_L}{l(l_L+jl)}=\frac{1}{l(l_L+jl)}\ge\frac{1}{l^2(j+1)}\ge\frac1{l^2(D+1)}=\frac{1}{l^2s(r)},
\]
where we used $l_L\le l$ and $j+1\le D+1$.
\end{proof}

\begin{proposition}\label{prop:sharp}
Let $x>0$ be rational and $u=k/l\ge0$ rational in lowest terms. Then for every $\gamma\in\Gamma$ with $\gamma(x)\ne u$,
\[
|\gamma(x)-u|\ge\frac{1}{l^2\,(s(x)+|\gamma|)}.
\]
In particular, the exponent $Q(u)=1$ in Theorem~\ref{thm:mobius} cannot be improved when $x$ and $u$ are rational.
\end{proposition}

\begin{proof}
Combine Lemma~\ref{lem:SB} with Lemma~\ref{lem:s}.
\end{proof}

\begin{remark}
We expect the same lower bound, with a constant depending on $x$, for irrational starting points $x$, but we have not proved it.
\end{remark}

\subsection{A general upper limit on the rate}
The following observation applies to any pair of maps.

\begin{proposition}\label{prop:count}
Let $T_0,T_1:(0,\infty)\to(0,\infty)$ be arbitrary, let $x>0$ and let $\mu>2$. Then for Lebesgue almost every $u\ge0$ there are only finitely many $\gamma\in\Gamma(T_0,T_1)$ with $|\gamma(x)-u|<\mu^{-|\gamma|}$.
\end{proposition}

\begin{proof}
Let $I_\gamma=(\gamma(x)-\mu^{-|\gamma|},\gamma(x)+\mu^{-|\gamma|})$. There are $2^n$ words of length $n$, so $\sum_\gamma|I_\gamma|=\sum_n2^n\cdot2\mu^{-n}<\infty$, and the Borel--Cantelli lemma shows that almost every $u$ lies in only finitely many $I_\gamma$.
\end{proof}

\section{The affine pair}\label{sec:affine}

Throughout this section $T_0(x)=x/a$ and $T_1(x)=bx+1$ with $a,b>1$, and $\sigma,\tau$ act as $T_0,T_1$. We put $\alpha=\log_ab>0$.

\subsection{A normal form}

\begin{lemma}\label{lem:normal}
Let $\gamma\in\Gamma$ contain $N$ letters $\tau$ and $M$ letters $\sigma$. For $0\le j<N$ let $\ell_j$ be the number of letters $\sigma$ to the left of the $(j+1)$-st letter $\tau$, counted from the left. Then $0\le\ell_0\le\ell_1\le\cdots\le\ell_{N-1}\le M$ and
\begin{equation}\label{eq:normal}
\gamma(x)=\frac{b^N}{a^M}\,x+\sum_{j=0}^{N-1}\frac{b^j}{a^{\ell_j}}\qquad(x\in\R).
\end{equation}
Conversely, for all integers $N\ge0$ and $0\le\ell_0\le\cdots\le\ell_{N-1}\le M$, the word
\[
\gamma=\sigma^{\ell_0}\tau\,\sigma^{\ell_1-\ell_0}\tau\cdots\sigma^{\ell_{N-1}-\ell_{N-2}}\tau\,\sigma^{M-\ell_{N-1}}
\]
has length $N+M$ and satisfies \eqref{eq:normal}.
\end{lemma}

\begin{proof}
Write $\gamma=w_1\cdots w_L$. Each $T_i$ is affine, $T_0'=1/a$ and $T_1'=b$, so $\gamma'=b^Na^{-M}$. The constant term of $\gamma$ is the sum, over the positions $i$ with $w_i=\tau$, of the constant $1$ of $T_1$ multiplied by the derivative of $w_1\cdots w_{i-1}$. If $w_i$ is the $(j+1)$-st $\tau$ from the left, this derivative is $b^ja^{-\ell_j}$. The converse is immediate.
\end{proof}

\subsection{The Kronecker-type estimate}
We use the following inhomogeneous form of Kronecker's theorem (\cite[Ch.~III, Thm.~II]{Cassels}, \cite[\S2.4]{Steuding}): if $\alpha$ is irrational and $\beta\in\R$, then there are infinitely many pairs of integers $(m,n)$ with $m>0$ and
\begin{equation}\label{eq:kron}
|m\alpha-n-\beta|<\frac3m.
\end{equation}

\begin{theorem}\label{thm:kronecker}
Suppose $\alpha=\log_ab$ is irrational, and put $t=a^{12}b^4$. Then for every $x>0$ and $u\ge0$ there are infinitely many $\gamma\in\Gamma$ with
\begin{equation}\label{eq:kron-thm}
|\gamma(x)-u|<(u+1)\big(t^{1/|\gamma|}-1\big).
\end{equation}
\end{theorem}

\begin{proof}
If $u=0$, then $\gamma=\sigma^m$ gives $\gamma(x)=xa^{-m}<t^{1/m}-1$ for $m$ large. Let $u>0$ and $\beta=\log_a\big((b-1)u\big)$. By \eqref{eq:kron} there are infinitely many pairs of positive integers $(m,n)$ with
\begin{equation}\label{eq:pf1}
|m\alpha-n-\beta|<\frac3m,
\end{equation}
and since $|\log_a(b^m-1)-m\alpha|=|\log_a(1-b^{-m})|<1/m$ for $m$ large,
\begin{equation}\label{eq:pf2}
\big|\log_a(b^m-1)-n-\beta\big|<\frac4m.
\end{equation}
Let $\eta=\sigma^n\tau^{m-1}$. Since $\tau^{m-1}(1)=(b^m-1)/(b-1)$, we have $\eta(1)=(b^m-1)/((b-1)a^n)$, and \eqref{eq:pf2} says $|\log_a\eta(1)-\log_a u|<4/m$. Hence
\[
|\eta(1)-u|<u\big(a^{4/m}-1\big).
\]
Now let $\gamma=\eta\tau\sigma^m=\sigma^n\tau^m\sigma^m$, so that $\gamma(x)=\eta(bxa^{-m}+1)$. Since $\eta$ is affine with derivative $b^{m-1}a^{-n}$,
\[
|\gamma(x)-\eta(1)|=\frac{b^{m}x}{a^{m+n}}\le (b-1)u\,x\,a^{3/m-m},
\]
using $b^ma^{-n}=a^{m\alpha-n}<(b-1)ua^{3/m}$ from \eqref{eq:pf1}. The right side is exponentially small in $m$, so for $m$ large
\[
|\gamma(x)-u|\le(b-1)uxa^{3/m-m}+u(a^{4/m}-1)<(u+1)(a^{4/m}-1).
\]
Finally, \eqref{eq:pf1} gives $n<(\alpha+1)m$ for $m$ large, hence $|\gamma|=2m+n<(\alpha+3)m$, and
\[
a^{4/m}<a^{4(\alpha+3)/|\gamma|}=\big(a^{12}b^{4}\big)^{1/|\gamma|}.\qedhere
\]
\end{proof}

\begin{remark}
Since $t^{1/n}-1\sim(\log t)/n$, Theorem~\ref{thm:kronecker} says that every orbit is $O(n^{-1})$-dense, with the explicit constant $|\gamma(x)-u|\lesssim 4(u+1)(3\log a+\log b)/|\gamma|$.
\end{remark}

\subsection{Rational $\log_ab$}

\begin{lemma}\label{lem:digits}
Let $a>1$ and $u>0$. There are non-negative integers $d_0,d_1,\dots$ with $d_0<u$ and $d_i\le\lceil a\rceil-1<a$ for $i\ge1$, such that
\begin{equation}\label{eq:digits}
0<u-\sum_{i=0}^kd_ia^{-i}\le a^{-k}\qquad\text{for all }k\ge0,
\end{equation}
and $d_k>0$ for infinitely many $k$.
\end{lemma}

\begin{proof}
Let $d_0=\lceil u-1\rceil$ and $\rho_0=u-d_0\in(0,1]$. Given $\rho_k\in(0,a^{-k}]$, let $d_{k+1}=\lceil\rho_ka^{k+1}-1\rceil$ and $\rho_{k+1}=\rho_k-d_{k+1}a^{-k-1}$. Since $0<\rho_ka^{k+1}\le a$, we get $0\le d_{k+1}<\rho_k a^{k+1}\le a$ and $\rho_{k+1}=a^{-k-1}(\rho_ka^{k+1}-d_{k+1})\in(0,a^{-k-1}]$. If $d_k=0$ for all $k\ge K$, then $\rho_k=\rho_K>0$ for $k\ge K$, contradicting $\rho_k\le a^{-k}\to0$.
\end{proof}

\begin{lemma}\label{lem:construct}
Suppose $a=b>1$ and let $d_0,\dots,d_k$ be non-negative integers with $d_k>0$. Then there is $\eta\in\Gamma$ with $\eta(1)=\sum_{i=0}^kd_ia^{-i}$, $|\eta|\le k+2\sum_{i=0}^kd_i$, and $\eta'=a^{-k}$.
\end{lemma}

\begin{proof}
Note $\tau\sigma(z)=z+1$. For $k=0$ take $\eta=(\tau\sigma)^{d_0-1}$. For the inductive step write $P=\sum_{i=0}^{k+1}d_ia^{-i}=d_0+Q/a$ with $Q=\sum_{i=0}^{k}d_{i+1}a^{-i}$; if $\mu(1)=Q$ with $|\mu|\le k+2\sum_{i\ge1}d_i$ and $\mu'=a^{-k}$, then $\eta=(\tau\sigma)^{d_0}\sigma\mu$ satisfies $\eta(1)=P$, $|\eta|\le k+1+2\sum_{i\ge0}d_i$ and $\eta'=a^{-k-1}$.
\end{proof}

\begin{theorem}\label{thm:rational}
Suppose $\log_ab=n_0/m_0$ with $m_0,n_0$ positive integers, and put $A=a^{m_0}=b^{n_0}$ and
\[
\rho=A^{1/((2A+1)\max\{m_0,n_0\})}>1.
\]
Then for every $x>0$ and $u\ge0$ there are $C=C(a,b,x,u)$ and infinitely many $\gamma\in\Gamma$ with
\[
|\gamma(x)-u|<C\rho^{-|\gamma|}.
\]
In particular, if $a=b$ one may take $\rho=a^{1/(2a+1)}$.
\end{theorem}

\begin{proof}
\emph{The case $a=b$.} If $u=0$, use $\gamma=\sigma^M$. Let $u>0$ and let $d_i$ be as in Lemma~\ref{lem:digits}. For each $k$ with $d_k>0$, Lemma~\ref{lem:construct} gives $\gamma_k$ with $\gamma_k(1)=\sum_{i\le k}d_ia^{-i}$, $\gamma_k'=a^{-k}$ and
\[
|\gamma_k|\le k+2d_0+2k(\lceil a\rceil-1)\le(2a+1)k+2u.
\]
Then
\[
|\gamma_k(x)-u|\le\gamma_k'|x-1|+|\gamma_k(1)-u|\le(|x-1|+1)a^{-k}\le(x+2)a^{2u/(2a+1)}\,a^{-|\gamma_k|/(2a+1)}.
\]
The values $\gamma_k(1)$ are strictly increasing in $k$, so the $\gamma_k$ are distinct.

\emph{The general case.} Let $c=\sum_{i=0}^{n_0-1}b^i$ and $\theta(z)=cz$. Then
\[
S_0:=\theta^{-1}\circ T_0^{m_0}\circ\theta:\ z\mapsto z/A,\qquad S_1:=\theta^{-1}\circ T_1^{n_0}\circ\theta:\ z\mapsto Az+1,
\]
because $T_1^{n_0}(y)=b^{n_0}y+c$. Applying the first case to $(S_0,S_1)$, with base $A$, to $v=x/c$ and $w=u/c$, we get infinitely many words $\mu$ in $S_0,S_1$ with $|\mu(v)-w|<C'A^{-|\mu|'/(2A+1)}$, where $|\mu|'$ is the length in $S_0,S_1$. Each such $\mu$ equals $\theta^{-1}\gamma\theta$ for a word $\gamma$ in $T_0,T_1$ with $|\gamma|\le\max\{m_0,n_0\}|\mu|'$, and $|\gamma(x)-u|=c\,|\mu(v)-w|$.
\end{proof}

\begin{remark}
The original version of this theorem gave the bound $\exp(-C|\gamma|/(u+1))$, whose exponent deteriorates as $u$ grows; the estimate $|\gamma_k|\le(2a+1)k+2u$ shows that $u$ only affects the constant. (Also, in the case $u=0$ the original choice of constant $C=1/2$ fails when $a<e^{1/2}$.)
\end{remark}

\subsection{Exponential density for all $a,b>1$}

\begin{theorem}\label{thm:exp}
Let $a,b>1$ and put
\[
\kappa=\log_a\frac{ab}{a-1},\qquad \lambda=\Big(\frac{a}{a-1}\Big)^{1/(1+\kappa)}=\exp\Big(\frac{\log a\,\log\frac{a}{a-1}}{\log\frac{a^2b}{a-1}}\Big)>1.
\]
Then for every $x>0$ and $u\ge0$ there are $C=C(a,b,x,u)$ and infinitely many $\gamma\in\Gamma$ such that
\begin{equation}\label{eq:exp}
|\gamma(x)-u|<C\lambda^{-|\gamma|}.
\end{equation}
In particular, every orbit is $O(\lambda^{-n})$-dense in $[0,\infty)$.
\end{theorem}

\begin{proof}
\emph{The case $u=0$.} We have $\log\lambda\le\log a$, because $\log\frac{a}{a-1}\le\log\frac{a^2b}{a-1}$. So $\gamma=\sigma^M$ gives $\gamma(x)=xa^{-M}\le x\lambda^{-M}$.

\emph{A greedy expansion.} Let $u>0$. Set $\ell_{-1}=0$ and $r_{-1}=u$. For $j=0,1,2,\dots$, given an integer $\ell_{j-1}$ and $r_{j-1}>0$, let $\ell_j^*$ be the least integer $\ell$ with $b^ja^{-\ell}<r_{j-1}$, and put
\[
\ell_j=\max\{\ell_{j-1},\ell_j^*\},\qquad t_j=b^ja^{-\ell_j},\qquad r_j=r_{j-1}-t_j.
\]
Since $\ell_j\ge\ell_j^*$ we have $t_j<r_{j-1}$, so $r_j>0$ and the recursion continues; also $0=\ell_{-1}\le\ell_0\le\ell_1\le\cdots$ and $r_{j}=u-\sum_{i\le j}t_i$. Call step $j$ \emph{free} if $\ell_j=\ell_j^*$.

(i) \emph{If step $j$ is free, then $t_j\ge r_{j-1}/a$; hence $r_j\le(1-1/a)r_{j-1}$ and $r_j\le(a-1)t_j$.} Indeed, by minimality of $\ell_j^*$, $a\,t_j=b^ja^{-(\ell_j^*-1)}\ge r_{j-1}$.

(ii) \emph{If step $j$ is free, so is step $j+1$.} We must show $\ell_{j+1}^*\ge\ell_j$, i.e.\ that $\ell=\ell_j-1$ violates $b^{j+1}a^{-\ell}<r_j$. But $b^{j+1}a^{-(\ell_j-1)}=ab\,t_j>(a-1)t_j\ge r_j$ by (i).

(iii) \emph{Some step is free.} If steps $0,\dots,j$ are not free, then $\ell_0=\cdots=\ell_j=0$ and $t_i=b^i$, so $1+b+\cdots+b^j<u$. Hence the first free step $J$ satisfies $b^J\le 1+(b-1)u$.

By (i)--(iii), for every $N>J$,
\begin{equation}\label{eq:rbounds}
r_{N-1}\le u\Big(1-\frac1a\Big)^{N-J}\qquad\text{and}\qquad r_{N-1}\le(a-1)t_{N-1}=\frac{a-1}{b}\,b^Na^{-\ell_{N-1}}.
\end{equation}

\emph{The words.} For $N>J$ let $M_N=\max\{\ell_{N-1},\lceil\kappa N\rceil\}$ and let $\gamma_N$ be the word of Lemma~\ref{lem:normal} with data $N$, $\ell_0\le\cdots\le\ell_{N-1}$ and $M=M_N$. Then $L:=|\gamma_N|=N+M_N$ and
\[
\gamma_N(x)=b^Na^{-M_N}x+u-r_{N-1},\qquad\text{so}\qquad|\gamma_N(x)-u|\le x\,b^Na^{-M_N}+r_{N-1}.
\]
Put $\vartheta=\kappa-\alpha=\log_a\frac{a}{a-1}>0$, so that $a^{-\vartheta}=1-1/a$ and $\lambda=a^{\vartheta/(1+\kappa)}$.

\emph{Case $M_N=\ell_{N-1}$.} By \eqref{eq:rbounds}, $|\gamma_N(x)-u|\le(x+\frac{a-1}b)\,a^{N\alpha-M_N}$. Since $M_N\ge\kappa N$ and $\vartheta+\alpha=\kappa$, a direct computation gives
\[
(1+\kappa)(M_N-N\alpha)-\vartheta(N+M_N)=(1+\alpha)(M_N-\kappa N)\ge0,
\]
that is, $M_N-N\alpha\ge\vartheta L/(1+\kappa)$. Hence $|\gamma_N(x)-u|\le(x+\frac{a-1}{b})\lambda^{-L}$.

\emph{Case $M_N=\lceil\kappa N\rceil>\ell_{N-1}$.} Then $b^Na^{-M_N}\le a^{-\vartheta N}=(1-1/a)^N$, and by \eqref{eq:rbounds},
\[
|\gamma_N(x)-u|\le\Big(x+u\big(\tfrac{a}{a-1}\big)^{J}\Big)\Big(1-\frac1a\Big)^N.
\]
As $L\le(1+\kappa)N+1$, we have $(1-1/a)^N=\lambda^{-(1+\kappa)N}\le\lambda^{1-L}$.

In both cases $|\gamma_N(x)-u|\le C\lambda^{-|\gamma_N|}$ with $C=\lambda\big(x+u(\frac a{a-1})^J\big)+\frac{a-1}b$, and the words $\gamma_N$, $N>J$, are distinct since they contain different numbers of letters $\tau$. Replacing $C$ by $2C$ gives strict inequality.
\end{proof}

\begin{remark}\label{rem:exp}
\begin{enumerate}[label=(\arabic*)]
\item Theorem~\ref{thm:exp} makes no use of the arithmetic of $\log_ab$, and it supersedes Theorem~\ref{thm:kronecker} and, for small $a$, Theorem~\ref{thm:rational}. For instance, $\lambda(2,3)\approx1.2133$ and $\lambda(2,2)=2^{1/3}$, whereas Theorem~\ref{thm:rational} gives $\rho=2^{1/5}$ when $a=b=2$; for $a=b=10$, on the other hand, Theorem~\ref{thm:rational} gives the larger base $10^{1/21}\approx1.116$ against $\lambda\approx1.053$.
\item By Proposition~\ref{prop:count}, for almost every $u$ the base $\lambda$ in \eqref{eq:exp} can never exceed $2$; thus exponential density is the optimal order, although we have not tried to optimize $\lambda$. Numerical experiments suggest that the greedy words of the proof perform considerably better than the bound $\lambda^{-|\gamma|}$.
\item The constant $C$ depends on $u$ only through $u(\frac a{a-1})^J\le u\,(1+(b-1)u)^{\log_b\frac{a}{a-1}}$, a polynomial in $u$.
\end{enumerate}
\end{remark}

\section{Prime approximation of real numbers}\label{sec:primes}

Let $p_n$ denote the $n$-th prime. Goldston, Pintz and Y\i ld\i r\i m \cite{GPY} (see also \cite{GMPY}) proved that
\begin{equation}\label{eq:gpy}
\liminf_{n\to\infty}\frac{p_{n+1}-p_n}{\log p_n}=0,
\end{equation}
equivalently $|p_{n+1}/p_n-1|=o(\log p_n/p_n)$ for infinitely many $n$. Zhang \cite{Zhang} proved that $\liminf(p_{n+1}-p_n)<\infty$, and Maynard \cite{Maynard} and the Polymath project \cite{Polymath} improved this to $\liminf(p_{n+1}-p_n)\le246$. Hence $|p_{n+1}/p_n-1|\le246/p_n$ for infinitely many $n$. This suggests the following.

\begin{conjecture}[Prime approximation]\label{conj:prime}
For every $u\ge0$ and $\eps>0$ there are infinitely many pairs of distinct primes $(p,q)$ with
\begin{equation}\label{eq:pac}
\abs{\frac pq-u}<\eps\,\frac{\log q}{q}.
\end{equation}
\end{conjecture}

The cases $u=0$ (take $p=2$) and $u=1$ (by the bounded gaps theorem) are known, the latter in the stronger form $|p/q-1|\le246/q$. A naive probabilistic model suggests that even $|p/q-u|<c/q$ should have infinitely many solutions for every $c>0$, since for each prime $q$ the interval $(uq-c,uq+c)$ should contain a prime with probability about $2c/\log q$, and $\sum_q1/\log q$ diverges.

Unconditionally, Baker, Harman and Pintz \cite{BHP} proved that $[x-x^{0.525},x]$ contains a prime for all sufficiently large $x$. Choosing a prime $p$ in $[uq-(uq)^{0.525},uq]$ shows that, for every fixed $u>0$ and \emph{every} sufficiently large $q$ (prime or not),
\begin{equation}\label{eq:bhp}
0\le u-\frac pq\le u^{0.525}q^{-0.475}.
\end{equation}
Assuming the Riemann Hypothesis, Cram\'er \cite{Cramer} showed $p_{n+1}-p_n=O(\sqrt{p_n}\log p_n)$, which improves \eqref{eq:bhp} to $|p/q-u|\le C_uq^{-1/2}\log q$, still far from \eqref{eq:pac}.

\subsection{Ratios of logarithms of primes}
We need a Mills-type lemma \cite{Mills,Tee}; we include the short proof.

\begin{lemma}\label{lem:mills}
Suppose $0<\beta<1$ and $p_{n+1}-p_n=O(p_n^\beta)$. Then for every real $u>1/(1-\beta)$ there is $\theta>1$ such that $\lfloor\theta^{u^n}\rfloor$ is prime for every $n\ge1$.
\end{lemma}

\begin{proof}
The hypothesis gives a prime in $[y,y+Ky^\beta]$ for all $y\ge2$ and some $K$. For a large integer $P$, the interval $[P^u,(P+1)^u-2]$ has length at least $uP^{u-1}-2$, which exceeds $K(P^u)^\beta$ for $P$ large because $u-1>u\beta$. So we may choose primes $P_1<P_2<\cdots$, with $P_1$ large, such that $P_n^u\le P_{n+1}\le(P_n+1)^u-2$. Then $P_n^{u^{-n}}$ is non-decreasing and $(P_n+1)^{u^{-n}}$ is strictly decreasing, so $\theta=\lim P_n^{u^{-n}}$ satisfies $P_n\le\theta^{u^n}<P_n+1$, i.e.\ $\lfloor\theta^{u^n}\rfloor=P_n$, for all $n\ge1$.
\end{proof}

\begin{theorem}\label{thm:logprimes}
For every $u>40/19$ there is an infinite sequence of primes $q_1<q_2<\cdots$ such that, for all $n$,
\begin{equation}\label{eq:logchain}
-\frac{1}{q_{n+1}\log q_n}<\frac{\log q_{n+1}}{\log q_n}-u<\frac{u}{q_n\log q_n}.
\end{equation}
Assuming the Riemann Hypothesis, the same holds for every $u>2$.
\end{theorem}

\begin{proof}
By \cite{BHP}, the hypothesis of Lemma~\ref{lem:mills} holds with $\beta=0.525$, and $1/(1-0.525)=40/19$. (Under the Riemann Hypothesis, Cram\'er's bound gives it for every $\beta>1/2$.) Let $q_n=\lfloor\theta^{u^n}\rfloor$. From $q_n\le\theta^{u^n}<q_n+1$ and $\theta^{u^{n+1}}=(\theta^{u^n})^u$ we get
\[
\frac{\log q_{n+1}}{\log(q_n+1)}<u<\frac{\log(q_{n+1}+1)}{\log q_n},\qquad q_{n+1}<(q_n+1)^u .
\]
The right inequality gives $u-\frac{\log q_{n+1}}{\log q_n}<\frac{\log(1+1/q_{n+1})}{\log q_n}<\frac1{q_{n+1}\log q_n}$. The left one gives
\[
\frac{\log q_{n+1}}{\log q_n}-u<\log q_{n+1}\,\frac{\log(1+1/q_n)}{\log q_n\log(q_n+1)}<\frac{u}{q_n\log q_n},
\]
using $\log q_{n+1}<u\log(q_n+1)$ and $\log(1+1/q_n)<1/q_n$.
\end{proof}

The feature of Theorem~\ref{thm:logprimes} is that the approximations form a \emph{chain}: each $q_n$ serves as the numerator of one approximation and the denominator of the next. For a single pair, a much stronger estimate holds for all $u>0$:

\begin{proposition}\label{prop:logpair}
Let $u>0$. For every sufficiently large prime $q$ there is a prime $p$ with
\[
0\le u-\frac{\log p}{\log q}\le\frac{2\,q^{-0.475u}}{\log q}.
\]
\end{proposition}

\begin{proof}
Let $X=q^u$ and pick a prime $p\in[X-X^{0.525},X]$ \cite{BHP}. Then $0\le\log X-\log p\le-\log(1-X^{-0.475})\le2X^{-0.475}$; divide by $\log q$.
\end{proof}

For $u>40/19$ we have $q^{-0.475u}=o(q^{-1})$, so Proposition~\ref{prop:logpair} is stronger than \eqref{eq:logchain} pair by pair; its proof, however, does not produce chains.

Finally, if $0<v<19/40$, applying Theorem~\ref{thm:logprimes} to $u=1/v$ gives a sequence of primes with
\[
\abs{\frac{\log q_n}{\log q_{n+1}}-v}<\frac{1}{q_n\log q_{n+1}}.
\]
This implies $q_n>q_{n+1}^v\exp(-1/q_n)$, and therefore, for every fixed $\eps>0$ and $n$ large,
\[
\abs{\frac{\log q_n}{\log q_{n+1}}-v}<\frac{1+\eps}{q_{n+1}^v\log q_{n+1}}.
\]

\end{document}